\documentclass[11pt,reqno]{amsart}
\usepackage{amssymb,graphicx}

\usepackage{amsfonts,amscd,amsthm,amsmath}

\numberwithin{equation}{section}

\begin{document}
\begin{center}
\textbf{\large{Closed-form formulae for the derivatives\\ of trigonometric functions at
rational multiples of $\pi$}}

\vskip 2mm\textbf{Djurdje Cvijovi\'{c}} \vskip 2mm
{\it Atomic
Physics Laboratory, Vin\v{c}a Institute of
Nuclear Sciences \\
P.O. Box $522,$ $11001$ Belgrade$,$ Republic of Serbia}\\

\textbf{E-Mail: djurdje@vinca.rs}\\

\vskip 2mm \begin{quotation}
\textbf{Abstract.} In this sequel to our recent note \cite{Cvijovic0}  it is shown, in a unified manner, by making use of some basic properties of certain special functions, such as the Hurwitz zeta function, Lerch zeta
function and Legendre chi function, $\;$that the values of all
 derivatives of four trigonometric functions at rational
multiples of $\pi$ can be expressed in closed form as simple finite
sums involving the Bernoulli and Euler polynomials. In addition, some particular
cases are considered.
\end{quotation}
\end{center}



\vskip 2mm\section{Introduction} Recently, Adamchik  \cite[p. 4, Eq. 26]{Adamchik} solved completely a long-standing problem of finding a closed-form
expression for the higher derivatives of the cotangent function
\cite{Adamchik,Apostol,Kolbig}. In a recent note \cite{Cvijovic0},
we have deduced  closed-form expressions for the values of the
higher derivatives of the cotangent function at rational multiples
of $\pi$, which are considerably simpler than those
found by K\"{o}lbig \cite[Theorem 4]{Kolbig}.

In this sequel, by using different arguments, we provide another
more direct and shorter proof  of this result (Theorem $1(i)$)
and also derive corresponding (and previously unknown) simple closed-form formulae for the cosecant, tangent and secant functions (Theorems $1(ii)$ and $2$). More specifically, we show in a unified manner (by
making use of the properties of certain special functions, such as
the Hurwitz zeta function, Lerch zeta function and Legendre chi
function) that the values of all derivatives of these trigonometric
functions at rational multiples of $\pi$ can be expressed as finite
sums involving the Bernoulli and Euler polynomials.

\vskip 2mm\section{Statement and proof of the results}

In what follows, we set $\imath:=\sqrt{-1}$, an empty sum  is
interpreted as zero, and the classical Bernoulli and
Euler polynomials, $B_n(x)$ and $E_n(x)$, are defined by \cite[pp.
59 and 63]{Srivastava}:
\begin{equation}
\aligned \frac{t\, e^{t x}}{e^t-1}=\sum_{n=0}^{\infty}
B_n(x)\frac{t^n}{n!}&\quad (\lvert t\rvert <2\pi)\quad
\textup{and}\quad \frac{2\, e^{t x}}{e^t+1}= \sum_{n=0}^{\infty}
E_n(x)\frac{t^n}{n!}\quad (\lvert t\rvert <\pi)\\
&(n\in \mathbb{N}_{0}:=\mathbb{N}\cup{0};\mathbb{N}:={1, 2, 3,
\ldots}).
\endaligned
\end{equation}

Our main results are as follows.

\medskip
\noindent{\bf Theorem 1.} {\em If}  $n$,  $p$,
 $q\in \mathbb{N}$ {\em and} $p$ {\em and} $q$ {\em are such that}
$1\leq  p<q,$ {\em then, in terms of the Bernoulli polynomials}
$B_n(x),$ {\em we have:}

\begin{align*}
&\textup(i)\quad
\left.\frac{\textup{d}^{n}}{\textup{d}x^{n}}\cot(\pi x) \,\right|_{
x=\frac{p}{q}} = \imath^{n+1} \,\frac{2\,(2 \pi q)^{n}}{n+1}\sum_{
\alpha\,=1}^q e^{\frac{2 \pi \imath  \alpha p}{q}} B_{ n+1} \left(
\tfrac{\alpha}{q}\right);
\\
&\textup(ii)\quad \left.\frac{\textup{d}^{n}}{\textup{d}
x^{n}}\csc(\pi x) \,\right|_{x=\frac{p}{q}}= \imath^{n+1}
\,\frac{2\,(2 \pi q)^{n}}{n+1}\sum_{ \alpha\,=1}^q e^{\frac{\pi
\imath (2  \alpha -1) p}{q}} B_{n+1}\left(\tfrac{2  \alpha-1}{2
\,q}\right).
\end{align*}
\vskip 2mm \noindent{\bf Theorem 2.} {\em If}  $n, p\in \mathbb{N},$ $q\in\mathbb{N}\setminus\{1,2\}$ {\em and} $p$ {\em and}
$q$ {\em are such that} $1\leq  p < q/2,$ {\em then, in terms of the
Bernoulli and Euler polynomials,} $B_n(x)$ {\em and} $E_n(x),$ {\em
we have:}
\begin{align*}
&\textup(i)\,\,\left.\frac{\textup{d}^{n}}{\textup{d} x^{n}}\tan(\pi
x)\right|_{x=\frac{p}{q}}=\imath^{n+1}\,(2 \pi q)^{n}\sum_{
\alpha\,=1}^q (-1)^{\alpha-1} e^{\frac{2 \pi \imath \alpha  p}{q}}
\left\{ \begin{array}{l} -
 E_n \left(\frac{ \alpha}{q}\right)\quad\quad\,\,\,\, q\,\textup{is odd} \\
 \frac{2}{n+1}\, B_{n+1} \left(\frac{ \alpha}{q}\right)\,\,\, q\,\textup{is even}
 \end{array} \right.;
\\
&\textup(ii)\,\,\left.\frac{\textup{d}^{n}}{\textup{d}
x^{n}}\sec(\pi x)\right|_{x=\frac{p}{q}} = \imath^{n} (2 \pi q)^{n}
\\
&\qquad\qquad\qquad\qquad\qquad\; \cdot \sum_{ \alpha\,=1}^q
(-1)^{\alpha-1} e^{\frac{\pi \imath (2
 \alpha-1) p}{q}} \left\{
\begin{array}{l}
 E_n \left(\frac{ 2 \alpha-1}{2\,q}\right)\qquad\qquad\; q\,\textup{is odd} \\
 - \frac{2}{n+1}\, B_{n+1} \left(\frac{ 2 \alpha-1}{2\,q}\right)\,\,\, q\,\textup{is even}
 \end{array} \right..
\end{align*}
\vskip 2mm \noindent{\bf Remark 1.} As the following
examples show, the above formulae can be rewritten in a
somewhat different form:
\begin{align*}
& \left.\csc(\pi x)^{(2 n -1)} \,\right|_{ x=\frac{p}{q}}=
\frac{(-1)^n (2 \pi q)^{2 n-1}}{n}\sum_{\alpha\,=1}^q B_{2
n}\left(\tfrac{2 \alpha-1}{2\,q}\right)\cos\left(\tfrac{\pi (2
\alpha-1)
 p}{q}\right),
\\
& \left.\csc(\pi x)^{(2 n)} \,\right|_{
x=\frac{p}{q}}=\frac{(-1)^{n-1} \,2(2 \pi q)^{2 n}}{2
n+1}\sum_{\alpha\,=1}^q B_{2n+1}\left(\tfrac{2
\alpha-1}{2\,q}\right)\sin\left(\tfrac{\pi (2 \alpha-1)
 p}{q}\right).
\end{align*}
\noindent Further, it should be noted that a number of simpler expressions could be obtained
by specializing the results for
$x=\frac{1}{2},\frac{1}{3},\frac{2}{3},\frac{1}{4},\frac{3}{4},\frac{1}{6}$
and $\frac{5}{6}$. For instance:
\begin{align*}
&\left.\csc(\pi x)^{(2 n-1)} \,\right|_{ x=\frac{1}{2}} = 0,\quad
\left.\csc(\pi x)^{(2 n)} \,\right|_{ x=\frac{1}{2}} = (-1)^{n-1}
\frac{4 (4 \pi)^{2 n}}{2 n+1} B_{2 n+1}\left(\tfrac{1}{4}\right);
\\
&\left.\cot(\pi x)^{(2 n-1)} \,\right|_{ x=\frac{1}{4}} =
\left.\cot(\pi x)^{(2 n-1)} \,\right|_{ x=\frac{3}{4}}= (-1)^n (4
\pi)^{2 n-1} (2^{2 n}-1) \frac{B_{2 n} (0)}{n};
\\
&\left.\cot(\pi x)^{(2 n)} \,\right|_{x=\frac{1}{4}} = -
\left.\cot(\pi x)^{(2 n)} \,\right|_{x=\frac{3}{4}}= (-1)^{n-1}
\frac{4 (8 \pi)^{2 n}}{2 n +1} B_{2 n+1}\left(\tfrac{1}{4}\right).
\end{align*}
\vskip 2mm Let $\zeta (s,a)$ and $\zeta^{*} (s,a)$ denote,
respectively, the Hurwitz (or generalized) zeta
function defined by \cite[p. 96]{Srivastava}:
\begin{equation}
\zeta(s,a):=\sum_{k\,=0}^\infty \frac{1}{(k + a)^s}\quad(a \notin
\mathbb{Z}_{0}^{-}:=\{0,-1, -2, -3, \ldots\};\,\Re{(s)}>1)
\end{equation}
\noindent and its alternating counterpart defined
by (see e.g. \cite[p. 761]{Chang}):
\begin{equation} \zeta^{*}(s,a):=\sum\limits_{k\, =
0}^{\infty }\frac{(-1)^{k}}{(k+a)^{s}}\qquad (a \notin
\mathbb{Z}_{0}^{-};\,\Re{(s)}>0),
\end{equation}
\noindent \noindent and observe that the following relationships hold \cite[p. 85, Eq. (17)]{Srivastava}:
\begin{equation}
\zeta(1-n,a)=-\frac{B_n(a)}{n}\qquad(n\in \mathbb{N})
\end{equation}
\noindent and  \cite[p. 761, Eq. (2.3)]{Chang}
\begin{equation}
\zeta^{*}(-n,a)=\frac{1}{2}\, E_n(a)\qquad(n\in \mathbb{N}_{0}).
\end{equation}

\vskip 2 mm \noindent {\bf Proof of Theorem  1.}  Our proof requires
the use of the Lerch zeta function and the Legendre chi
function which are defined, respectively, by the series \cite[p. 89, Eq.
(7)]{Srivastava}:
\begin{equation}\ell_s(\xi):=\sum_{k\,=1}^{\infty} \frac{e^{2 \pi
\imath k \xi}}{k^{s}}\qquad (\xi\in \mathbb{R};\Re{(s)}>1),
\end{equation}
\noindent and (see, for instance, \cite{Cvijovic1}):
\begin{equation}
\chi_s(z):=\sum_{k\,=0}^{\infty}\frac{z^{2 k+1}}{(2 k+1)^s}
\qquad(\lvert z\rvert\leq 1; \Re{(s)}>1),
\end{equation}
\noindent and their meromorphic continuations over the whole
$s$-plane.
First, we shall show that
\begin{equation}\aligned
\ell_0(\xi)=-\frac{1}{2}+ \frac{\imath}{2}\cot(\pi\xi)&,\qquad
\ell_{1-n} (\xi)=\frac{\imath}{2 (2\pi\imath)^{n-1}}
\frac{\textup{d}^{n-1}
}{\textup{d}\xi^{n-1}}\,\cot(\pi\xi)\\
&\hspace{-5mm}(\xi\in\mathbb{R}\setminus
\mathbb{Z};\,n\in\mathbb{N}\setminus\{1\})\endaligned
\end{equation}
\noindent and
\begin{equation}\chi_{1-n} (e^{\pi\imath\, \xi})=\frac{\imath}{2 (\pi\imath)^{n-1}}  \frac{\textup{d}^{n-1}
}{\textup{d}\xi^{n-1}}\,\csc(\pi\xi)\qquad(\xi\in\mathbb{R}\setminus
\mathbb{Z};\,n\in\mathbb{N}).
\end{equation}
To prove (2.8) note that
\begin{equation}\frac{\partial }{\partial \xi}\, \ell_s(\xi)=2\pi\imath
\,\ell_{s-1}(\xi),
\end{equation}
\noindent which, in turn, follows from  (2.6) for $\Re{(s)}>2$ and
by analytic continuation for all $s$. The definition  in (2.6) also
yields $\ell_1(\xi)=-\log\left(1-e^{2\pi\imath\,\xi}\right) $
 $(\xi\in\mathbb{R}\setminus \mathbb{Z})$ and from
this we obtain $\ell_0(\xi)$  by (2.10). Using (2.10) repeatedly with
initial value $\ell_0(\xi)$ leads to $\ell_{1-n}(\xi)$ given by
(2.8). Likewise, we have (2.9) by making use of
\begin{equation}\frac{\partial }{\partial \xi}\, \chi_s(e^{\pi\imath\, \xi})=\pi\imath
\,\chi_{s-1}(e^{\pi\imath\, \xi})
\end{equation}
\noindent and
\begin{equation}
\chi_0(e^{\pi\imath\, \xi})=\frac{\imath}{2}\csc(\pi\xi)
\qquad(\xi\in\mathbb{R}\setminus \mathbb{Z}).
\end{equation}

Second, from (2.6), we obtain ({\em cf.} \cite[p. 1531]{Cvijovic2})
\begin{equation*}
\ell_s\left(\tfrac{p}{q}\right)= \sum_{k = \,0}^{\infty} \frac{e^{2
\pi \imath \,(k + 1) p/q}}{(k + 1)^s} = \sum_{\alpha =
\,0}^{q-1}\sum_{k =\,0}^{\infty} \frac{e^{2 \pi \imath k p} \,e^{2
\pi \imath (\alpha +1 ) p/q}}{q^s (k + (\alpha+1)/q)^s},
\end{equation*}
\noindent so that, in view of the definition in (2.2), we have (see
\cite[p. 1530, Eq. 8(a)]{Cvijovic2})
\begin{equation}
\ell_s\left(\tfrac{p}{q}\right)= \frac{1}{q^s}\sum_{\alpha\, =1}^q
\zeta\left(s,\tfrac{\alpha}{q}\right)\,e^{\frac{2 \pi \imath \alpha
p}{q}}\qquad(p,  q\in\mathbb{N}; p=1,\ldots, q).
\end{equation}
\noindent Similarly, starting from (2.7), we find that
\begin{equation}
\chi_s\left(e^{\frac{\pi\imath p}{q}}\right) = \frac{1}{(2
q)^s}\sum_{\alpha\, =1}^q \zeta\left(s,\tfrac{2 \alpha-1}{2
q}\right)\,e^{\frac{\pi \imath (2 \alpha - 1) p}{q}}\qquad(p,
q\in\mathbb{N}; p =1,\ldots, q).
\end{equation}
\noindent It should be noted that (2.13) and (2.14) are derived for
$\Re{(s)}>1$ but, by the principle of analytic continuation,  they
hold true for any complex $s$, $s\neq 1$.

Lastly, set $s =1-n$ $(n\in\mathbb{N}\setminus\{1\})$. Part
$(i)$ now follows at once by applying (2.13) in conjunction with
(2.8) and (2.4).  Also, Equation (2.14), in
conjunction with (2.9) and (2.4),  gives part $(ii)$. \qed

\vskip 2 mm \noindent {\bf Proof of Theorem  2.} Instead of the functions $\ell_s(\xi)$ and $\chi_s(z)$  used above, we now introduce the functions  $\ell^{*}_s(\xi)$ and $\chi^{*}_s(z)$ by means of the series
\begin{equation}\ell^{*}_s(\xi):=\sum_{k\,=1}^{\infty} (-1)^{k-1}\,\frac{e^{2 \pi
\imath k \xi}}{k^{s}}\qquad (\xi\in \mathbb{R};\Re{(s)}>0)
\end{equation}
\noindent and
\begin{equation}
\chi^{*}_s(z):=\sum_{k\,=0}^{\infty}(-1)^k\,\frac{z^{2 k+1}}{(2
k+1)^s} \qquad(\lvert z\rvert\leq 1; \Re{(s)}>0),
\end{equation}
\noindent and their meromorphic continuations over the whole
$s$-plane.

The proof follows precisely along the same lines as that of
Theorem 1: We first show that
\begin{equation}\aligned
\ell^{*}_0(\xi)=& \frac{1}{2}+ \frac{\imath}{2}\tan(\pi\xi),\qquad
\ell^{*}_{1-n} (\xi)=\frac{\imath}{2 (2\pi\imath)^{n-1}}
\frac{\textup{d}^{n-1}
}{\textup{d}\xi^{n-1}}\,\tan(\pi\xi)\\
&\quad\left(\xi\in\mathbb{R}\setminus \{(2 k + 1)\tfrac{1}{2}\,\vert
k\in\mathbb{Z}\};\,n\in\mathbb{N}\setminus\{1\}\right)\endaligned
\end{equation}
\noindent and
\begin{equation}\chi^{*}_{1-n}(e^{\pi\imath\, \xi})=\frac{\imath}{2 (\pi\imath)^{n-1}}\frac{\textup{d}^{n-1}
}{\textup{d}\xi^{n-1}}\,\sec(\pi\xi)\quad(\xi\in\mathbb{R}\setminus
\{(2 k + 1)\tfrac{1}{2}\,\vert k\in\mathbb{Z}\};\,n\in\mathbb{N}).
\end{equation}
\noindent We next establish the following formulae valid for $p,
q\in\mathbb{N}$ $(p =1,\ldots, q)$:
\begin{equation}
\ell^{*}_s\left(\tfrac{p}{q}\right)=
\frac{1}{q^s}\sum_{\alpha\,=1}^q (-1)^{\alpha-1}e^{\frac{2 \pi
\imath \alpha p}{q}} \,\left\{
\begin{array}{l} \zeta^{*}\left(s,\frac{\alpha}{q}\right)\qquad q\,\textup{is odd} \\
 \zeta\left(s,\frac{\alpha}{q}\right)\qquad\;\, q\,\textup{is even}
 \end{array} \right.,
\end{equation}
\begin{equation}
\chi^{*}_s\left(\tfrac{p}{q}\right)=
\frac{1}{q^s}\sum_{\alpha\,=1}^q (-1)^{\alpha -1}e^{\frac{\pi \imath
(2 \alpha - 1) p}{q}} \,\left\{
\begin{array}{l} \zeta^{*}\left(s,\frac{2 \alpha -1}{2 q}\right)\qquad q\,\textup{is odd} \\
 \zeta\left(s,\frac{2 \alpha -1}{2 q}\right)\qquad\;\, q\,\textup{is even}
 \end{array} \right.
\end{equation}
\noindent and finally set $s = 1-n$ $(n\in\mathbb{N}
\setminus\{1\})$. It is straightforward to see that part $(i)$
is obtained by (2.19) together with (2.15) and either (2.4) or
(2.5). Also, Equation (2.20), together with (2.16) and either
(2.4) or (2.5), yields part $(ii)$.\qed

\end{document}